\documentclass[a4paper,english,12pt]{amsart}
\title[Countable eq.~relations and quotient spaces]{Countable Borel equivalence relations and quotient Borel spaces}

\usepackage[english]{babel}
\usepackage[latin1]{inputenc}
\usepackage{latexsym}
\usepackage{amsmath,amsthm}
\usepackage{amsfonts,dsfont,amssymb,mathrsfs,bbm,bbold}
\usepackage[all]{xy}

\swapnumbers
\newtheorem{theo}{Theorem}[section]
\newtheorem{prop}[theo]{Proposition}

\newtheorem{lemma}[theo]{Lemma}
\newtheorem{dfnt}[theo]{Definition}

\theoremstyle{remark}
\newtheorem{remark}[theo]{Remark}
\newtheorem{claim}{Claim}[theo]
\newtheorem{ex}[theo]{Example}

\numberwithin{equation}{section}

\newcommand{\cal}[1]{\ensuremath{\mathcal{#1}}}
\renewcommand{\bf}[1]{\ensuremath{\mathbf{#1}}}
\renewcommand{\frak}[1]{\ensuremath{\mathfrak{#1}}}

\newcommand{\R}{\ensuremath{\mathds{R}}}

\newcommand{\N}{\ensuremath{\mathds{N}}}
\newcommand{\Z}{\ensuremath{\mathds{Z}}}

\newcommand{\res}[2]{\ensuremath{{#1}\!\upharpoonright\!{#2}}}
\newcommand{\set}[1]{\ensuremath{\left\{#1\right\}}}
\newcommand{\seqn}[1]{\ensuremath{\left(#1\right)}}
\newcommand{\abs}[1]{\ensuremath{\left\lvert#1\right\rvert}}

\newcommand{\expset}[2]{\ensuremath{{}^{#1}#2}}
\newcommand{\embeds}{\ensuremath{\sqsubseteq}}

\newcommand{\BSigma}[2]{\ensuremath{\boldsymbol{\Sigma}^{#1}_{#2}}}
\newcommand{\BPi}[2]{\ensuremath{\boldsymbol{\Pi}^{#1}_{#2}}}

\newcommand{\Bor}{\ensuremath{\mathbf{Bor}}}

\DeclareMathOperator{\dom}{dom}
\DeclareMathOperator{\rng}{rng}
\DeclareMathOperator{\proj}{proj}
\DeclareMathOperator{\id}{id}

\DeclareMathOperator{\End}{End}
\DeclareMathOperator{\Aut}{Aut}
\DeclareMathOperator{\Simm}{\frak{S}}

\newenvironment{bigtheocite}[2]{\refstepcounter{theo}\begin{trivlist}\item{\normalfont\thetheo}{\bfseries.\ #1\ \textnormal{(#2)}.}\itshape}{\end{trivlist}}

\newenvironment{tfae}{\begin{list}{\roman{tfaecounter}.}{\setlength{\labelwidth}{\leftmargin}}\usecounter{tfaecounter}}{\end{list}}

\providecommand{\multichotomycase}[1]{\setcounter{multycounter}{#1}\setcounter{multysubcounter}{0}}

\providecommand{\romanval}[1]{\setcounter{tempcounter}{#1}\roman{tempcounter}}

\providecommand{\implication}[2]{(\textit{\romanval{#1}$\,\rightarrow\,$\romanval{#2}})}

\providecommand{\parref}[1]{\mbox{\S\ref{#1}}}

\providecommand{\commentout}[1]{}

\newcounter{tempcounter}
\newcounter{tfaecounter}
\newcounter{multycounter}\newcounter{multysubcounter}

\begin{document}

\author[R.~Pinciroli]{Roberto Pinciroli}
\address{Scuola Normale Superiore\\56100 Pisa\\Italy}
\email{r.pinciroli@sns.it}
\thanks{The results in this paper emerged from various interesting discussions I had with prof. A.~Louveau, to whom I express my sincere gratitude.}

\begin{abstract}
We consider countable Borel equivalence relations on quotient Borel spaces. We prove a generalization of the Feldman\textendash Moore representation theorem, but provide some examples showing that other very simple properties of countable equivalence relations on standard Borel spaces may fail in the context of nonsmooth quotients.
\end{abstract}

\maketitle


\section{Quotient Borel spaces}

We recall some basic definitions: a \emph{standard Borel space} is a measurable space $(X,\Bor(X))$ which is isomorphic to the Borel space of a complete separable metric space; a \emph{Borel equivalence relation} on $X$ is an equivalence relation $E$ which is Borel as a subset of $X\times X$, and $E$ is \emph{countable} (\emph{finite}) if all its equivalence classes are countable (finite).
Given two equivalence relations $E,F$ on two sets $X,Y$ respectively, a \emph{morphism} $f$ of $E$ to $F$ (denoted as $f:E\preceq F$) is a map $f:X\to Y$ sending $E$-equivalent points to $F$-equivalent points; $f$ is a \emph{reduction} if moreover $xEy\ \leftrightarrow\ f(x)Ff(y)$: in this case we write $f:E\leq F$. If we consider Borel morphisms and reductions, we write $f:E\preceq_B F$ and $f:E\leq_B F$; a \emph{monomorphism} is an injective morphism (denoted as $f:E\subseteq F$) and an \emph{embedding} is an injective reduction ($f:E\embeds F$).

A Borel equivalence relation $E$ on a standard Borel space $X$ is \emph{smooth} if it is Borel reducible to the equality relation on some standard Borel space $Y$. For a countable Borel eq.~relation, the following properties are equivalent (see \cite{JKLCountableBorelEqRelations}, \cite{KMTopicsOrbitEquivalence}):
\begin{tfae}
\item $E$ is smooth;
\item $E$ admits a countable Borel separating family;
\item there is a Borel transversal $T\subseteq X$ for $E$;
\item there is a Borel selector $\varphi:X\to X$ for $E$;
\end{tfae}
more generally, for an arbitrary Borel $E$ one has \romanval{1}$\leftrightarrow$\romanval{2} and \romanval{3}$\leftrightarrow$\romanval{4}, but the latter ones are stronger than the others (see \cite{HKLGlimmEffrosDichotomy}; in order to restore the inverse implication, one may replace ``Borel'' with ``\bf{C}-measurable'' in \romanval{3} and \romanval{4}).
Another characterization of smoothness is given by the kind of Borel structure which a countable equivalence relation induces on the quotient set:
\begin{lemma}\label{lemma:QuotientBeingStandard}
Let $E$ be a countable Borel equivalence relation on a standard Borel space $X$. The following statements are equivalent:
\begin{tfae}
\item $E$ is smooth;
\item the quotient $\sigma$-algebra $\Bor(X/E)$ is countably generated;
\item $(X/E,\Bor(X/E))$ is a standard Borel space.
\end{tfae}
\end{lemma}
\proof
\implication{3}{2} By definition: just consider a countable base for a separable metrizable topology inducing the Borel structure of $X/E$.
\implication{2}{1} If \set{A_n\ :\ n\in\N} is a countable generating family for $\Bor(X/E)$ and $\pi:X\to X/E$ is the quotient projection, the set of preimages \set{\pi^{-1}[A_n]} is a countable Borel separating family for $E$.
\implication{1}{3} Let $f$ be a reduction $(X,E)\leq_B (Y,=)$, where $Y$ is some standard Borel space; since $E$ is countable, $f$ is a countable-to-1 Borel function between Polish spaces, so it maps Borel sets to Borel sets. It follows that the induced function $\tilde{f}:X/E\to Y$ is a Borel isomorphism of $X/E$ onto the Borel subset $f[X]\subseteq Y$, equipped by the relative $\sigma$-algebra $\res{\Bor(Y)}{f[X]}$, which is a standard Borel space.
\qed

By the preceding lemma, one has to be very cautious when dealing with quotients by \emph{nonsmooth} equivalence relations, since the usual theorems about standard (or just separable metrizable) Borel spaces are likely to fail in this more general context. To keep as much information as possible about the original ``nice'' structure, our official definition will be as follows:

\begin{dfnt}[Quotient Borel space]\label{def:QuotientBorelSpace}
A \emph{quotient Borel space} is a couple $(X,E)$, where $X$ is a standard Borel space and $E$ is a countable Borel equivalence relation on $X$. The underlying set of $(X,E)$ is just the quotient set $X/E$.
\end{dfnt}

A Borel subset of $(X,E)$ is simply an element of the quotient $\sigma$-algebra $\Bor\big(\frac{X}{E}\big)$; in other words, letting $\pi_E:X\to X/E$ be the quotient projection, $A\subseteq(X,E)$ is Borel if and only if $\pi_E^{-1}[A]$ is a Borel subset of $X$. Similarly, $A\subseteq(X,E)$ is \BSigma{1}{1} iff $\pi_E^{-1}[A]$ is \BSigma{1}{1} in $X$, and so on for the other projective classes.

However, for $n\geq 1$, it turns out that the ``right'' notion of Borelness for $(n+1)$-ary relations and for $n$-ary functions between quotient Borel spaces is \emph{not} the usual one from the context of measurable spaces and maps: here again we want to exploit the original standard Borel structures.

\begin{dfnt}[Products]\label{def:ProductsQuotientSpaces}
Let $(X,E)$ and $(Y,F)$ be quotient Borel spaces as in definition \ref{def:QuotientBorelSpace}.
The \emph{product} $(X,E)\times(Y,F)$ is the quotient Borel space $(X\times Y,E\times F)$, where
\begin{equation}\label{eq:ProductEquivalenceRelations}
(x_1,y_1)\ E\times F\ (x_2,y_2) \qquad\text{iff}\qquad (x_1\ E\ x_2\;\;\text{and}\;\; y_1\ F\ y_2);
\end{equation}
the underlying set of $(X,E)\times(Y,F)$ is naturally identified with the product of the respective underlying sets, $\frac{X}{E}\times\frac{Y}{F}$.
\end{dfnt}

This definition straightforwardly generalizes to that of product of any finite number of quotient Borel spaces.
\begin{remark}\label{remark:ProductQuotients}
It is easy to verify that the $\sigma$-algebra $\Bor\big(\frac{X\times Y}{E\times F}\big)$ contains the product $\sigma$-algebra $\Bor\big(\frac{X}{E}\big)\otimes\Bor\big(\frac{Y}{F}\big)$; in fact, we shall show (proposition \ref{prop:ProductNonsmoothQuotients}) that, whenever one of the two spaces $(X,E)$ and $(Y,F)$ is nonsmooth and the other is uncountable, the inclusion is always \emph{strict}.
\end{remark}

\begin{dfnt}[Relations and functions]\label{def:RelationsFunctionsQuotientSpaces}
A Borel $n$-ary relation $R$ on $(X_1,E_1)$~\ldots\ $(X_n,E_n)$ is a Borel subset $R\subseteq(X_1,E_1)\times\ldots\times(X_n,E_n)$.
A Borel $n$-ary function $f:(X_1,E_1)\times\ldots\times(X_n,E_n)\to(Y,F)$ is a function $\frac{X_1}{E_1}\times\ldots\times\frac{X_n}{E_n}\to\frac{Y}{F}$ which is Borel as an $(n+1)$-ary relation $f\subseteq(X_1,E_1)\times\ldots\times(X_n,E_n)\times(Y,F)$.
\end{dfnt}

These definitions allow for some basic properties of Borel relations and functions on standard Borel spaces to be preserved also for quotient spaces (see lemmata \ref{lemma:PropertiesBorelFunctions} and \ref{lemma:CategoryQuotientSpaces}); we also have a general ``lifting'' property which can be applied to avoid dealing with quotients (lemma \ref{lemma:ExistenceLiftings}).
The proofs make essential use of the countability assumption on the equivalence relations.

\begin{lemma}\label{lemma:PropertiesBorelFunctions}
Let $R\subseteq (X,E)\times(Y,F)$ be a Borel relation and $f:(X,E)\to (Y,F)$ be a Borel function between quotient Borel spaces.
\begin{enumerate}
\item if $R$ has countable sections (i.e. for all $x\in (X,E)$ the sections $R_x=\set{y\ :\ (x,y)\in R}$ are countable), then $\dom R$ is Borel in $(X,E)$;
\item the inverse image $f^{-1}[B]$ of a Borel subset $B\subseteq(Y,F)$ is Borel in $(X,E)$;
\item if $f$ is countable-to-1, the image $f[A]$ of a Borel subset $A\subseteq(X,E)$ is Borel in $(Y,F)$;
\end{enumerate}
\end{lemma}
\proof
(1.)
Let $\tilde{R}$ be the Borel subset $\pi_{E\times F}^{-1}[R]\subseteq X\times Y$ and note that $\dom\tilde{R}=\pi_E^{-1}[\dom R]$; since $F$ is countable and $R$ has countable sections, $\tilde{R}$ has countable sections too, hence $\dom\tilde{R}$ is Borel (see theorem \ref{bigtheo:LusinNovikov}). By the definition of Borel subset of the quotient space, $\dom R$ is Borel.

(2.)
If $B$ is Borel in $(Y,F)$ then $\frac{X}{E}\times B$ is Borel in $(X,E)\times(Y,F)$; it follows that $\big(f\cap\frac{X}{E}\times B\big)$ is Borel and $f^{-1}[B]=\dom\big(f\cap\frac{X}{E}\times B\big)$ is Borel by (1.) (here each section contains at most one point). The proof of (3.) is similar: just consider the Borel relation with countable sections $\big(f^{-1}\cap\frac{Y}{F}\times A\big)\subseteq(Y,F)\times(X,E)$ and apply (1.) to $f[A]=\dom\big(f^{-1}\cap\frac{Y}{F}\times A\big)$.
\qed

\begin{lemma}[Existence of liftings]\label{lemma:ExistenceLiftings}
Let $X,Y$ be standard Borel spaces, $E,F$ countable Borel equivalence relations on $X$ and $Y$ respectively. A function $f:(X,E)\to (Y,F)$ is Borel if and only if it has a Borel lifting $\tilde{f}:X\to Y$.
\end{lemma}
\begin{equation*}
\xymatrix{
X\ar[d]_{\pi_E}\ar@{-->}[r]^{\tilde{f}} &Y\ar[d]^{\pi_F}\\
(X,E)\ar[r]_{f} &(Y,F)
}
\end{equation*}
\proof
Suppose $f:(X,E)\to(Y,F)$ is Borel: by definition, the set $P=\pi_{E\times F}^{-1}[f]$ is Borel in $X\times Y$ and has countable sections, because $f$ is a function and $F$ is countable. Applying the Lusin\textendash Novikov theorem \ref{bigtheo:LusinNovikov}, let $\tilde{f}$ be a Borel uniformization of $P$: clearly $\tilde{f}$ is a lifting of $f$.

Conversely, assuming that $\tilde{f}$ is any lifting of $f$, $\pi_{E\times F}^{-1}[f]$ is just the $(E\times F)$-saturation of $\tilde{f}$ in $X\times Y$. By the Feldman\textendash Moore theorem \ref{bigtheo:FeldmanMooreRepresentationTheorem}, there exists a countable group $\varGamma$ acting in a Borel way on $X\times Y$ with orbit equivalence $E\times F$, therefore $\pi_{E\times F}^{-1}[f]=\bigcup_{\gamma\in\varGamma} \gamma\cdot\tilde{f}$ is Borel (hence $f$ is Borel) if $\tilde{f}$ is Borel.
\qed

\begin{lemma}\label{lemma:CategoryQuotientSpaces}
The composition of two Borel functions between quotient Borel spaces is Borel: the classes of quotient Borel spaces and of Borel functions form a concrete category.
Given objects $(X_1,E_1),\ldots,(X_n,E_n)$, the product $(X,E)=(X_1,E_1)\times\ldots\times(X_n,E_n)$, along with the projection maps $\proj_i:(X,E)\to(X_i,E_i)$, is in fact the categorical product: each projection is Borel and, for every space $(Y,F)$ and every family of morphisms $f_i:(Y,F)\to(X_i,E_i)$, the $n$-uple $f=(f_1,\ldots,f_n)$ is a morphism $(Y,F)\to(X,E)$.
\end{lemma}
\proof
These facts easily follow from the lifting lemma \ref{lemma:ExistenceLiftings}. As an example, we prove Borelness of composition: given Borel maps $f:(X,E)\to(Y,F)$ and $g:(Y,F)\to(Z,G)$, pick Borel liftings $\tilde{f}:X\to Y$ and $\tilde{g}:Y\to Z$; the composite $\tilde{g}\circ\tilde{f}$ is then a Borel function lifting $g\circ f$, therefore the latter is Borel.
\qed

\begin{prop}\label{prop:ProductNonsmoothQuotients}
Let $X,Y$ be standard Borel spaces, $E,F$ countable Borel equivalence relations on $X$ and $Y$ respectively; assume moreover that $F$ is nonsmooth and $X$ (or equivalently $X/E$) is uncountable. Then there is a Borel subset of $(X,E)\times(Y,F)$ which is not measurable as a subset of the product measurable space $\frac{X}{E}\times\frac{Y}{F}$.
\end{prop}
\proof
By the Glimm\textendash Effros dichotomy (see \cite{HKLGlimmEffrosDichotomy}) there is a Borel injection $g:(\expset{\omega}{2},E_0)\to(Y,F)$, where $E_0$ is the nonsmooth equivalence relation of eventual equality between $2$-valued sequences,
\begin{equation}\label{eq:RelationE0}
\seqn{x_n}\ E_0\ \seqn{y_n}\qquad\text{iff}\qquad\forall^{\infty} n\ (x_n=y_n)
\end{equation}
(here $\forall^{\infty}n$ means $\exists n_0\ \forall n\geq n_0$); by the Silver dichotomy (\cite{SilCountingNumberEquivalenceClasses}, \cite{MKInfiniteGamesEffectiveDST}), there exists also a Borel injection $f:\expset{\omega}{2}\to(X,E)$.
Let $\pi_0$ be the quotient map $\expset{\omega}{2}\to(\expset{\omega}{2},E_0)$;
the product $f\times g$ is an injective Borel morphism $\expset{\omega}{2}\times(\expset{\omega}{2},E_0)\to(X,E)\times(Y,F)$ (lemma \ref{lemma:CategoryQuotientSpaces}), so the image of the Borel relation $\pi_0$ is a Borel subset of $(X,E)\times(Y,F)$ (lemma \ref{lemma:PropertiesBorelFunctions}).
We claim that $(f\times g)[\pi_0]$ is not a measurable subset of $\frac{X}{E}\times\frac{Y}{F}$. Suppose otherwise: since both $f:\expset{\omega}{2}\to\frac{X}{E}$ and $g:\frac{\expset{\omega}{2}}{E_0}\to\frac{Y}{F}$ are measurable functions, we would have that $\pi_0$ is measurable in $\expset{\omega}{2}\times\frac{\expset{\omega}{2}}{E_0}$, hence $\pi_0$ would belong to some sub-$\sigma$-algebra generated by countably many ``vertical stripes'' \set{A_l\times\frac{\expset{\omega}{2}}{E_0}} and countably many ``horizontal'' ones \set{\expset{\omega}{2}\times B_m}, with $A_l$ Borel in \expset{\omega}{2} and $B_m$ Borel in $\frac{\expset{\omega}{2}}{E_0}$ for $l,m\in\N$. Let $x,y$ be two $E_0$-inequivalent points of \expset{\omega}{2}: since $(x,[x]_{E_0})\in\pi_0$ but $(x,[y]_{E_0})\notin\pi_0$, for some index $m$ we should have $[x]_{E_0}\in B_m\ \leftrightarrow\ [y]_{E_0}\notin B_m$: in other words, \set{\pi_0^{-1}[B_m]} is a countable Borel separating family for $E_0$, a contradiction.
\qed

There are, however, important facts about relations in standard Borel spaces which cease to be true for quotient Borel spaces; among these, one is of particular interest in connection with countable equivalence relations (see the discussion of enumerable relations and the Feldman\textendash Moore theorem in \parref{sec:CountableEqRel}):

\begin{bigtheocite}{The Lusin\textendash Novikov uniformization theorem}{see \cite{KecClassicalDST}}\label{bigtheo:LusinNovikov}
Let $X,Y$ be standard Borel spaces and $P\subseteq X\times Y$ be Borel with countable sections. Then $P$ has a Borel uniformization $\varphi_P$ and $\dom\varphi_P=\dom P=\proj_X[P]$ is Borel.
Moreover $P$ can be written as a countable union of Borel graphs.
\qed
\end{bigtheocite}

This uniformization property doesn't hold anymore in the context of quotient Borel spaces.
\begin{ex}\label{ex:NonexistenceUniformization}
Let $E$ be any nonsmooth countable Borel equivalence relation on a standard Borel space $X$ and $\pi_E:X\to(X,E)$ be the quotient projection. The inverse relation of the graph of $\pi_E$ is a Borel subset of $(X,E)\times X$ with countable sections (they are exactly the $E$-equivalence classes in $X$); nonetheless, for every uniformization  $\varphi:(X,E)\to X$, the composition $\varphi\circ\pi_E$ is a selector for $E$, so $\varphi$ cannot be Borel.
\end{ex}


\section{Countable Borel equivalence relations}
\label{sec:CountableEqRel}

We begin our study of countable eq. relations on quotient spaces isolating some particular subclasses, according to the cardinality of the equivalence classes, to the generation or representability properties and, finally, to the regularity.

\begin{dfnt}
We say that an equivalence relation $E$ on a set $X$ has \emph{index $\leq n$} if all its equivalence classes have cardinality less than or equal to $n$; $E$ is \emph{finite} (or has index $<\aleph_0$) if its equivalence classes are finite.
\end{dfnt}

\begin{dfnt}
An equivalence relation $F$ on a quotient Borel space $(X,E)$ is \emph{(Borel) countably generated} (Borel $\aleph_0$-generated) if there is a countable family \set{f_n\ :\ n\in\N} of Borel endomorphisms of $(X,E)$ such that $F$ is the smallest subset of $(X,E)^2=(X^2,E^2)$ being an equivalence relation and containing every $f_n$.
$F$ is \emph{(Borel) enumerable} if there is a countable family \set{f_n\ :\ n\in\N} of Borel endomorphisms of $(X,E)$ such that $F=\bigcup_nf_n$; \set{f_n} is then called an \emph{enumeration} of $F$.
\end{dfnt}

A (Borel) enumerable equivalence relation is clearly (Borel) $\aleph_0$-generated; if the family \set{f_n\ :\ n\in\N} is an enumeration of $F$, then each $f_n$ is a countable-to-1 function and, for all points $x$, \set{f_n(x)} is an enumeration of the equivalence class $[x]_F$ of $x$.
We remark that, even on a standard Borel space, an equivalence relation $F$ which is generated by countably many Borel functions may well be uncountable and that, in this case, $F$ might be a complete \BSigma{1}{1} set (example \ref{ex:CtblyGenAnalyticEqRel}); however, if $F$ is already countable, then every function $f\subseteq F$ is countable-to-1, and lemma \ref{lemma:BorelCtblyGeneratedAnalytic} below applies.

\begin{ex}\label{ex:TailEquivalenceEndomorphism}
Given a Borel endomorphism $f:X\to X$, the equivalence relation generated by $f$ is simply the \emph{tail equivalence relation} of $f$, whose classes are the \emph{grand-orbits} of $f$:
\begin{equation}\label{eq:TailEquivalenceEndomorphism}
E(f)\;=\;E_t(f)\;=\;\set{(x,y)\ :\ \exists m,n\in\N\ \big(f^m(x)=f^n(y)\big)}.
\end{equation}
In this case $E(f)$ is clearly Borel, and $E(f)$ is countable exactly when $f$ is countable-to-1.
\end{ex}

\begin{ex}\label{ex:CtblyGenAnalyticEqRel}
Fix a complete \BSigma{1}{1} subset $A\subseteq\expset{\omega}{\omega}$, a Borel function $f_0:\expset{\omega}{\omega}\to\expset{\omega}{\omega}$ with $\rng f_0=A$ and a point $x_0\in A$.
Pick $X=\expset{\omega}{\omega}\times 2$ and define two Borel endomorphisms $f,g$ of $X$ as follows:
\begin{equation*}
f(x,i)\;=\;%
\begin{cases}
(f_0(x),1) &\ \text{if}\ i=0,\\
(x,1) &\ \text{if}\ i=1;
\end{cases}
\qquad\qquad
g(x,i)\;=\;%
\begin{cases}
(x_0,1) &\ \text{if}\ i=0,\\
(x,1) &\ \text{if}\ i=1.
\end{cases}
\end{equation*}
Let $E=E(f,g)$ be the equivalence relation generated by $f$ and $g$ on $X$; its equivalence classes are just the singletons \set{(x,1)} for $x\notin A$ and the whole remaining subset $(\expset{\omega}{\omega}\times\set{0})\cup(A\times\set{1})$.
$E$ is clearly an analytic subset of $X\times X$; since the Borel function
\begin{equation*}
\rho\;:\;\expset{\omega}{\omega}\to X\times X\;:\;x\mapsto\big((x,1),(x_0,1)\big)
\end{equation*}
reduces $A$ to $E$, $E$ is a complete \BSigma{1}{1} set.
\end{ex}

\begin{lemma}\label{lemma:BorelCtblyGeneratedAnalytic}
Let $(X,E)$ be a quotient Borel space and \set{f_n\ :\ n\in\N} a countable family of Borel endomorphisms of $(X,E)$. The equivalence relation $F$ generated by \set{f_n} is \BSigma{1}{1}; if every $f_n$ is countable-to-1, then $F$ is countable Borel.
\end{lemma}
\proof
For each $r\in\N$, consider the following subset $F_r$ of $(X,E)^{r+1}$:
\begin{multline*}
F_r\;=\;\set{(x_0,\ldots,x_r)\ :\ \forall k<r\ \exists j\in\N \phantom{\big|} \right.\\
\left.\big(x_{k+1}=x_k\ \ \text{or}\ \ x_{k+1}=f_j(x_k)\ \ \text{or}\ \ x_k=f_j(x_{k+1})\big)};
\end{multline*}
by definition, each $F_r$ is Borel. Moreover, if every $f_n$ is countable-to-1, for each choice of a point $x=x_0\in X$ there are only countably many possibilities for $x_1$, and then for $x_2$ and so on, in order for $(x_0,\ldots,x_r)$ to be in $F_r$, hence each section $(F_r)_x$ is countable.
It is straightforward to check that the (analytic) sets
\begin{equation*}
F_{(r)}\;=\;\set{(x,y)\ :\ \exists(z_0,\ldots,z_r)\in F_r\ (z_0=x\ \text{and}\ z_r=y)}\;\subseteq\;(X,E)\times (X,E)
\end{equation*}
are symmetric and reflexive binary relations on $(X,E)$ and that their union is transitive, so the equivalence relation $F$ generated by \set{f_n} coincides with $\bigcup_r F_{(r)}$ and so is \BSigma{1}{1}. When the $f_n$'s are countable-to-1, by the above remarks every set $F_{(r)}$ is the projection of a Borel set with countable sections, and its sections $(F_{(r)})_x$ are countable too; it follows from lemma \ref{lemma:PropertiesBorelFunctions} that every $F_{(r)}$, and thus also $F$, are countable Borel.
\qed

\begin{dfnt}
Let $\varGamma$ be a countable group acting on a quotient Borel space $(X,E)$. We say that the action is Borel if, for every $\gamma\in\varGamma$, the induced permutation $(\gamma\cdot)$ of $X/E$ is a Borel automorphism of $(X,E)$; in this case $(X,E)$ is a \emph{Borel $\varGamma$-space}. The \emph{orbit equivalence relation} $E_{\varGamma}^{(X,E)}$ induced on the (Borel) $\varGamma$-space $(X,E)$ is just the equivalence whose classes are the orbits of the action:
\begin{equation*}
x E_{\varGamma}^{(X,E)} y\quad\text{if and only if}\quad\exists\gamma\in\varGamma\ (y=\gamma\cdot x).
\end{equation*}
\end{dfnt}


\begin{dfnt}
A countable Borel equivalence relation $F$ on a quotient Borel space $(X,E)$ is \emph{smooth} if it admits a Borel transversal $T\subseteq(X,E)$.
\end{dfnt}
As when the space $(X,E)$ is standard, $F$ is smooth if and only if it has a Borel selector: given a Borel transversal $T$, the set
\begin{equation*}
\varphi\;=\;\set{(x,y)\ :\ xFy\ \text{and}\ y\in T}
\end{equation*}
is a Borel selector and, conversely, if $\varphi'$ is a Borel selector for $F$, the set $T'=\set{x\ :\ x=\varphi'(x)}$ is easily a Borel transversal.
Note that the usual definition of smoothness, which requires the existence of a countable Borel separating family (or, equivalently, reducibility to equality on some standard Borel space), is no longer suitable for equivalence relations on arbitrary quotient Borel spaces: by lemma \ref{lemma:QuotientBeingStandard}, the very identity relation on $(X,E)$ fails to satisfy these properties unless $E$ itself is smooth.

\begin{figure}[h]
\begin{equation*}
\vcenter{
\xymatrix{
\save [].[rrr]!C*[F--]\frm{}="top" \restore
\txt{Countable}\ar@{<->}[r] &\txt{(Countable)\\$\aleph_0$-generated}\ar@{<->}[r]
&\txt{Enumerable}\ar@{<->}[r] &\txt{Orbit equivalence\\of a countable group}\\
&\save "top"+<0pt,-10ex>*\txt{(Countable) smooth}="center"\ar"top" \restore &\\
&\save "top"+<0pt,-20ex>*\txt{Finite}\ar"center" \restore&\\
}}
\end{equation*}
\caption{Classes of countable eq.~rel.~on standard Borel spaces}
\label{fig:ClassesCountableEqRelSmooth}
\end{figure}
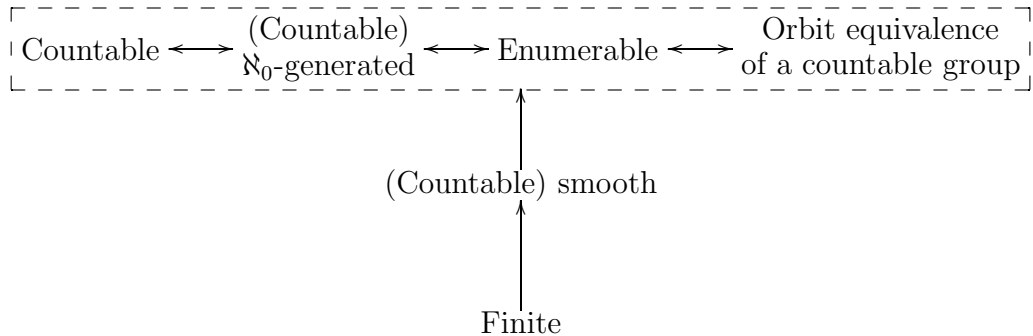

When the underlying space $(X,E)$ is smooth, the previous classes are very simply arranged as in figure \ref{fig:ClassesCountableEqRelSmooth}.

First of all, the entire collection of countable Borel equivalence relations gets collapsed to the class of orbit equivalences induced by Borel actions of countable groups: this is just the main statement of the classical Feldman\textendash Moore theorem, whose proof follows.
\begin{bigtheocite}{The Feldman\textendash Moore representation theorem}{see also \cite{KMTopicsOrbitEquivalence}}\label{bigtheo:FeldmanMooreRepresentationTheorem}
Let $E$ be a countable Borel equivalence relation on a standard Borel space $X$. Then there is a countable (discrete) group $\varGamma$ acting in a Borel way on $X$ with orbit equivalence $E_{\varGamma}^X=E$. The group $\varGamma$ and the action can be chosen so that
\begin{equation*}
xEy\;\leftrightarrow\;\exists \gamma\in\varGamma\ (\gamma^2=1\ \wedge\ \gamma\cdot x=y).
\end{equation*}
\end{bigtheocite}
\proof(Louveau)
The two key ingredients of the classical proof are the following facts:
\begin{enumerate}
\item any countable Borel equivalence relation $E$ on a standard Borel space $X$ is Borel enumerable, by the Lusin\textendash Novikov uniformization theorem \ref{bigtheo:LusinNovikov}: we can write $E=\bigcup_{m<\omega} f_m$, with $f_m$ Borel endomorphisms of $X$;
\item the $\sigma$-algebra $\Bor(X)$ of the Borel subsets of $X$ is countably generated: let $\cal{A}=\set{A_p\ :\ p<\omega}$ be a countable generating subalgebra of $\Bor(X)$.
\end{enumerate}
For each triple $(m,n,p)\in\omega^3$ define
\begin{equation*}
f_{m,n,p}(x)\;=\;%
\begin{cases}
f_m(x) &\ \text{if}\ x\in A_p,\ f_m(x)\notin A_p,\ f_n\circ f_m(x)=x,\\
f_n(x) &\ \text{if}\ x\notin A_p,\ f_n(x)\in A_p,\ f_m\circ f_n(x)=x,\\
x &\ \text{otherwise};
\end{cases}
\end{equation*}
it is clear that every $f_{m,n,p}$ is a Borel involution contained in $E$. Moreover, given a pair of distinct $E$-equivalent elements $x,y\in X$, we can find $m,n<\omega$ such that $f_m(x)=y$ and $f_n(y)=x$ (since the $f_m$'s cover $E$), and a $p<\omega$ such that $x\in A_p$, $y\notin A_p$ (since $\Bor(X)$ separates points and $\cal{A}$ is a generating subalgebra), hence $f_{m,n,p}(x)=y$.
The countable subgroup $\varGamma\leq\Aut(X)$ generated by \set{f_{m,n,p}\ :\ m,n,p<\omega} has the desired properties.
\qed

Another nice fact is that all finite Borel equivalence relations on standard Borel spaces are smooth: by the Isomorphism Theorem, there is no loss of generality if we just work with a finite Borel equivalence relation $F$ on a Borel subset $X$ of \R: in this case the set of minima of all the $F$-classes,
\begin{equation*}
T\;=\;\set{x\in X\ :\ \forall y\in[x]_F (x\leq y)}\;=\;\set{x\in X\ :\ \forall\gamma\in\varGamma\ (x\leq\gamma\cdot x)},
\end{equation*}
where $\varGamma$ is a countable group acting on $X$ with $F=E^X_{\varGamma}$, as provided by the Feldman\textendash Moore theorem, is easily seen to be a Borel transversal for $F$.

\begin{figure}[h]
\begin{equation*}
\vcenter{
\xymatrix@C=2cm{
&\txt{Countable}&\\
\txt{Finite}\ar[ur] &\txt{(Countable)\\countably generated}\ar[u] &\\
\txt{Index $\leq 3$}\ar[u]\ar@{.>}^-{\txt{\small(ex.~\ref{ex:NotCtblyGeneratedIndex3EqRel})}}|-{\bigotimes}[ru] &\txt{Enumerable}\ar[u] &\txt{(Countable)\\smooth}\ar[ul]\ar@{.>}^-{\txt{\small(ex.~\ref{ex:SmoothNonEnumerableIndex3EqRel})}}|-{\bigotimes}[l]\\
&\txt{Orbit equivalence\\of a countable group}\ar@{<->}[u]\save[u].[]*[F--]\frm{}\restore &\\
\txt{Index $\leq 2$}\ar[uu]\ar[ur]\ar@{.>}`r/2pt[rru]_-{\txt{\small(ex.~\ref{ex:NonsmoothIndex2EqRel})}}|-{\bigotimes}[rruu]&&
}}
\end{equation*}
\caption{Classes of countable eq.~rel.~on quotient Borel spaces}
\label{fig:ClassesCountableEqRelQuotientSpaces}
\end{figure}
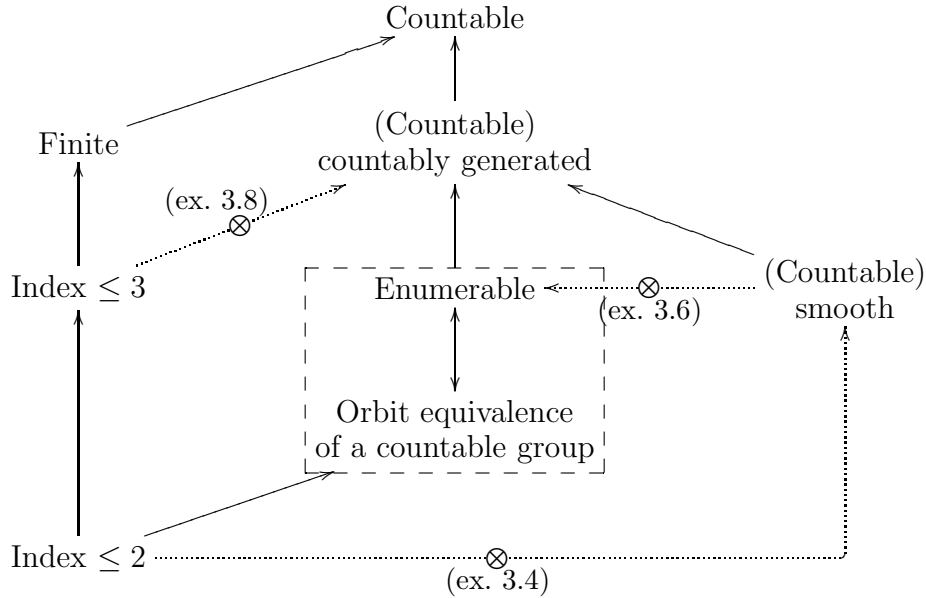

For general quotient Borel spaces, the previous facts are no longer valid and the overall picture is quite more complicated (see figure \ref{fig:ClassesCountableEqRelQuotientSpaces}). We just make a couple of simple observations: every smooth equivalence relation is Borel countably generated (in fact, $1$-generated, since one only needs a Borel selector to produce the entire equivalence), and each Borel eq. rel. $F$ of index $\leq\!\!2$ is induced by a single Borel involution,
\begin{equation*}
f(x)\;=\;%
\begin{cases}
y &\ \text{if}\ x\neq y\ \text{and}\ (x,y)\in F,\\
x &\ \text{if}\ x\notin\dom(F\smallsetminus\id);
\end{cases}
\end{equation*}
in particular, every index-$(\leq\!\!2)$ equivalence relation is an orbit equivalence $E_{\Z/2\Z}$ for some Borel action of $\Z/2\Z$.
The situation changes drastically already for index-$3$ relations (see examples \ref{ex:SmoothNonEnumerableIndex3EqRel} and \ref{ex:NotCtblyGeneratedIndex3EqRel}): it may happen that it is impossible to associate, in a ``simply definable'' and uniform way, to each point $x$ an element in the same equivalence class \emph{different from} $x$.

Since the Lusin\textendash Novikov theorem \ref{bigtheo:LusinNovikov} fails for quotient spaces, one cannot hope to prove a representation theorem like \ref{bigtheo:FeldmanMooreRepresentationTheorem} for all countable equivalence relations: for instance, example \ref{ex:SmoothNonEnumerableIndex3EqRel} provides a countable ($1$-generated) Borel equivalence relation which is not enumerable. Moreover, the usual proof of theorem \ref{bigtheo:FeldmanMooreRepresentationTheorem} makes use of the fact that the $\sigma$-algebra $\Bor(X)$ is countably generated when $X$ is standard, but this becomes false for nonsmooth quotients.
Nevertheless, we now proceed to show the remarkable result that the representation theorem is \emph{still} valid, for Borel enumerable relations, even without the assumption of the existence of a countable generating family for the Borel $\sigma$-algebra:

\begin{theo}[``Feldman\textendash Moore'' representation theorem for enumerable equivalence relations on quotient Borel spaces]\label{theo:FeldmanMooreQuotientSpaces}
Let $F$ be a Borel enumerable equivalence relation on a quotient Borel space $(X,E)$. Then there is a countable group $\varGamma$ acting in a Borel fashion on $(X,E)$ with orbit equivalence $E^{(X,E)}_{\varGamma}=F$.
\end{theo}

We shall use the following common notation: given an equivalence relation $F$ on a Borel space $X$, we write $[F]$ for the \emph{full group} of $F$, that is the subgroup of $\Aut X$ containing those automorphisms whose graph is a subset of $F$, and $[[F]]$ for the set of Borel injective graphs contained in $F$.

By the enumerability hypothesis, we can write $F$ as a countable union $\bigcup_{m<\omega}\varphi_m$ of Borel functions $\varphi_m\in\End(X,E)$. Fix a bijection $\omega\to\omega^2:n\mapsto((n)_0,(n)_1)$ and let, for each natural number $n$, $\psi_n=\varphi_{(n)_0}\cap\varphi_{(n)_1}^{-1}$, where $\varphi_m^{-1}$ denotes the inverse relation of $\varphi_m\subseteq(X/E)^2$: by construction, every $\psi_n$ is a Borel partial injection and
\begin{equation}\label{eq:FeldmanMooreQuotientSpaces1}
F\;=\;\bigcup_{n<\omega}\psi_n, \qquad\qquad \forall n<\omega\ \big(\psi_n\in[[F]]\big).
\end{equation}

The following lemma is the key to the proof:
\begin{lemma}\label{lemma:FeldmanMooreKeyLemma}
Any Borel partial injection $g_0\in[[F]]$ can be covered with two Borel automorphisms $g',g''\in[F]$.
\end{lemma}
\proof
Starting with $g_0$, we build an increasing sequence of partial injections $g_n$ by successively adjoining pieces of the various $\psi_n$: precisely, define by recurrence
\begin{equation}\label{eq:FeldmanMooreKeyLemma1}
g_{n+1}\;=\;g_n\cup\big(\psi_n\smallsetminus(\dom g_n\times (X/E))\smallsetminus ((X/E)\times \rng g_n)\big),
\end{equation}
that is, we attach to $g_n$ as much of $\psi_n$ as possible in order for $g_{n+1}$ to be again a partial injection. Put $g=\bigcup_{n<\omega} g_n$: $g$ is also in $[[F]]$.

\begin{claim}\label{claim:FeldmanMooreKeyLemma1}
If $y,z\in X/E$ and $yFz$ then either $y\in\dom g$ or $z\in\rng g$.
\end{claim}
Suppose not: by \eqref{eq:FeldmanMooreQuotientSpaces1} we can find an index $n$ such that $(y,z)\in\psi_n$, thus
\begin{align*}
(y,z)\;&\in\;\psi_n\smallsetminus(\dom g\times (X,E))\smallsetminus ((X,E)\times \rng g)\;\subseteq\\
&\subseteq\;\psi_n\smallsetminus(\dom g_n\times (X,E))\smallsetminus ((X,E)\times \rng g_n)\;\subseteq\;g_{n+1}\;\subseteq\;g,
\end{align*}
so $y\in\dom g$ and $z\in\rng g$, which is a contradiction.

Consider now the following family of Borel subsets of $X$, indexed by integers $n\in\Z$:
\begin{equation}\label{eq:FeldmanMooreKeyLemma2}
\begin{aligned}
&X_1\;=\;\dom g\smallsetminus\rng g, & &X_{n+1}\;=\;g[X_n];\\
&X_{-1}\;=\;\rng g\smallsetminus\dom g, & &X_{-(n+1)}\;=\;g^{-1}[X_{-n}];\\
&X_0\;=\;X\smallsetminus\bigcup_{n\neq 0} X_n.
\end{aligned}
\end{equation}
We plan to build $g'$ and $g''$ out of $g$ and $g^{-1}$ by means of a ``Schr\"oder\textendash Bern\-stein'' type argument, and the first step to accomplish this goal is the following

\begin{claim}\label{claim:FeldmanMooreKeyLemma2}
\seqn{X_n\ :\ n\in\Z} is a Borel partition of $(X,E)$.
\end{claim}
By the definition of $X_0$, one needs only check that $X_k$ and $X_l$ are disjoint for distinct nonzero $k,l\in\Z$; if $k$ and $l$ have the same sign, this is an easy consequence of the injectivity of $g$ and $g^{-1}$. Suppose then $k<0<l$ and $x\in X_k\cap X_l$: by \eqref{eq:FeldmanMooreKeyLemma2},
\begin{align*}
&z\;=\;g^{-(l-1)}(x)\;=\;g^{-l+1}(x)\;\in\;X_1, & &\text{so}\ z\notin\rng g,\\
&y\;=\;(g^{-1})^{-(\abs{k}-1)}(x)\;=\;g^{-k-1}(x)\;\in\;X_{-1}, & &\text{so}\ y\notin\dom g;
\end{align*}
however we have $yFxFz$ (since $g,g^{-1}\in[[F]]$) and therefore $yFz$ by transitivity, which contradicts claim \ref{claim:FeldmanMooreKeyLemma1}.

We proceed to construct the two desired bijections $g',g''\in[F]$ as follows,
\begin{align*}
g'\;=\;&\res{g}{X_0}\ \cup\ \bigcup_{n<\omega}\res{g}{X_{2n+1}}\ \cup\ \bigcup_{n<\omega}\res{g^{-1}}{X_{2n+2}}\ \cup\\
&\cup\ \bigcup_{n<\omega}\res{g^{-1}}{X_{-2n-1}}\ \cup\ \bigcup_{n<\omega}\res{g}{X_{-2n-2}},\\
g''\;=\;&\res{g^{-1}}{X_0}\ \cup\ \id_{X_1\cup X_{-1}}\ \cup\ \bigcup_{n<\omega}\res{g}{X_{2n+2}}\ \cup\ \bigcup_{n<\omega}\res{g^{-1}}{X_{2n+3}}\ \cup\\
&\cup\ \bigcup_{n<\omega}\res{g^{-1}}{X_{-2n-2}}\ \cup\ \bigcup_{n<\omega}\res{g}{X_{-2n-3}},
\end{align*}
(notice that $g$ and $g^{-1}$ are bijective on $X_0$) and we obtain
\begin{equation*}
g_0\cup g_0^{-1}\;\subseteq\;g\cup g^{-1}\;\subseteq\;g'\cup g''.
\qedhere
\end{equation*}

\proof[Proof of theorem \ref{theo:FeldmanMooreQuotientSpaces}]
By lemma \ref{lemma:FeldmanMooreKeyLemma}, every partial injection $\psi_n$ in \eqref{eq:FeldmanMooreQuotientSpaces1} is covered by two automorphisms $\psi'_n,\psi''_n\in[F]$, so the subgroup $\varGamma\leq[F]$ generated by all $\psi'_n$ and all $\psi''_n$ works.
\qed

\begin{remark}
The automorphisms $g'$, $g''$ given by lemma \ref{lemma:FeldmanMooreKeyLemma} may not be involutions of the space $(X,E)$, so we cannot deduce the last statement of the classical Feldman\textendash Moore theorem \ref{bigtheo:FeldmanMooreRepresentationTheorem}, that $F$ is covered by countably many involutions (induced by elements of $\varGamma$ of order $2$): in fact, this is false in general, see example \ref{ex:NotInvolutionGeneratedIndex3EqRel}. 
However, if both $g_0$ and all the partial injections $\psi_n$ in \eqref{eq:FeldmanMooreQuotientSpaces1} were involutions, then clearly the $g$ built in lemma \ref{lemma:FeldmanMooreKeyLemma} would be too; since $\dom g$ equals $\rng g$ in this case, by claim \ref{claim:FeldmanMooreKeyLemma1} $g$ would already be an automorphism of $(X,E)$. Hence the proof of theorem \ref{theo:FeldmanMooreQuotientSpaces} can be adapted to give a new proof of theorem \ref{bigtheo:FeldmanMooreRepresentationTheorem}.
\end{remark}

Let us conclude this section with a simple, weak uniformization lemma: here again the hypothesis of countability, which is sufficient for relations in standard Borel spaces, has to be strengthened (example \ref{ex:NonexistenceUniformization}):

\begin{lemma}[Weak uniformization]\label{lemma:WeakUniformization}
Let $R\subseteq(X,E)\times(Y,F)$ be a Borel binary relation between quotient Borel spaces. If $R$ is contained in an enumerable relation $S$ (i.e. $R$ is covered by countably many Borel functions $(X,E)\to(Y,F)$), then $R$ has a Borel uniformization.
\end{lemma}
\proof
Is straightforward: let \set{f_n\ :\ n\in\N} be a countable family of Borel functions, $f_n:(X,E)\to(Y,F)$, such that $R\subseteq\bigcup_nf_n$. Put
\begin{align*}
&\varPhi_R\;=\;\set{(x,n)\in(X,E)\times\N\ :\ (x,f_n(x))\in R\ \text{and}\ \forall k<n\ \big((x,f_k(x))\notin R\big)},\\
&\varphi_R\;=\;\set{(x,y)\in R\ :\ \exists n\in\N\ \big((x,n)\in\varPhi_R\ \text{and}\ y=f_n(x)\big)};
\end{align*}
clearly both $\varPhi_R$ and $\varphi_R$ are Borel functions, $\dom\varPhi_R=\dom\varphi_R=\dom R$ and $\varphi_R$ is a uniformization of $R$.
\qed


\section{Some examples; $2$-valued measures and free actions}
\label{sec:Examples}

We begin this section with a simple general criterion (proposition \ref{prop:CriterionFreeMeasurePreservingActions}) which is useful to build examples of nonsmooth and non-countably generated Borel equivalence relations.

\begin{lemma}
\label{lemma:ExistenceOfBorelCocyclesFreeActions}
Let $F=E_{\varGamma}^{(X,E)}$ be the orbit equivalence of some Borel \emph{free} action of a countable group $\varGamma$ on a quotient Borel space $(X,E)$. Then the action admits a Borel cocycle, i.e. a Borel function $\theta:F\to\varGamma$ satisfying the following properties:
\begin{enumerate}
\item for all $(x,y)\in F$, $y=\theta(x,y)\cdot x$;
\item for all triples of $F$-equivalent elements $xFyFz$, $\theta(x,z)=\theta(y,z)\theta(x,y)$.
\end{enumerate}
\end{lemma}
\proof
For each $\gamma\in\varGamma$, put $A_{\gamma}=\set{(x,y)\ :\ y=\gamma\cdot x}$: since $F=E_{\varGamma}$ and the action is Borel and free, \seqn{A_{\gamma}\ :\ \gamma\in\varGamma} is a Borel partition of $F$, so the required cocycle is simply
\begin{equation*}
\theta(x,y)\ =\ \gamma \qquad\text{iff}\qquad (x,y)\in A_{\gamma}.
\qedhere
\end{equation*}

\begin{prop}\label{prop:CriterionFreeMeasurePreservingActions}
Let $(X,E)$ be a nonsmooth quotient Borel space, $\mu$ a $2$-valued nonatomic Borel measure on $(X,E)$, $\varGamma$ a countable group acting in a Borel, free and $\mu$-preserving way on $(X,E)$, $\varDelta$ a subgroup of $\varGamma$.
\begin{enumerate}
\item If $\varGamma$ is not the trivial group, $E_{\varGamma}^{(X,E)}$ is nonsmooth;
\item if $F\subseteq E_{\varGamma}^{(X,E)}$ is a countably generated equivalence relation over $E_{\varDelta}^{(X,E)}$, then $F$ is $\mu$-almost contained in $E_{N_{\varGamma}(\varDelta)}^{(X,E)}$, where $N_{\varGamma}(\varDelta)$ is the normalizer of $\varDelta$ in $\varGamma$.
\end{enumerate}
\end{prop}
In the last assertion, the meaning of ``$\mu$-almost contained'' is that the set
$$\dom\big(F\smallsetminus E_{N_{\varGamma}(\varDelta)}^{(X,E)}\big)\;=\;\set{x\in(X,E)\ :\ [x]_F\nsubseteq [x]_{E_{N_{\varGamma}(\varDelta)}^{(X,E)}}}
$$(which is Borel by lemma \ref{lemma:PropertiesBorelFunctions}) is $\mu$-null.
\begin{remark}
The existence of a $2$-valued nonatomic Borel measure $\mu$ on $(X,E)$ is in fact equivalent to $E$ being nonsmooth: for one direction, if $E$ were smooth then $\Bor(X/E)$ would be countably generated (lemma \ref{lemma:QuotientBeingStandard}), so $\mu$ would take value $1$ on some atom of $\Bor(X/E)$, contradicting non-atomicity of $\mu$; for the other, recall that whenever $E$ is nonsmooth there is a Borel $E$-ergodic $E$-nonatomic measure $\tilde{\mu}$ on $X$, so it is sufficient to let $\mu=(\pi_E)_{*}\tilde{\mu}$, where $\pi_E$ denotes the quotient map $X\to X/E$.
This is just the countable case of the Harrington\textendash Kechris\textendash Louveau dichotomy for arbitrary Borel equivalence relations (see \cite{HKLGlimmEffrosDichotomy}), which was first proved by Effros \cite{EffTransformationGroupsCAlgebras}, \cite{EffPolishTransfGroupsClassification} and Weiss \cite{WeiMeasurableDynamics}.
\end{remark}
\proof
(1.)
Let $\theta:F\to\varGamma$ be a Borel cocycle associated to the action of $\varGamma$ on $(X,E)$ (lemma \ref{lemma:ExistenceOfBorelCocyclesFreeActions}). Suppose that $E_{\varGamma}^{(X,E)}$ is smooth over $E$ and let $\varphi:(X,E)\to(X,E)$ be a Borel selector; the function
\begin{equation}\label{eq:CocycleSelector}
\tilde{\varphi}\ :\ (X,E)\to\varGamma\ :\ x\mapsto\theta(x,\varphi(x))
\end{equation}
is Borel and satisfies, for all $x\in(X,E)$, $\varphi(x)=\tilde{\varphi}(x)\cdot x$. As $\gamma$ varies over $\varGamma$, the sets
\begin{equation*}
Z_{\gamma}\;=\;\tilde{\varphi}^{-1}(\gamma)\;=\;\set{x\in (X,E)\ :\ \tilde{\varphi}(x)=\gamma}
\end{equation*}
form a Borel partition of $(X,E)$ into countably many pieces.
If $\delta\in\varGamma$ is not the identity element (we're assuming $\varGamma$ nontrivial, hence such a $\delta$ exists), for every $\gamma\in\varGamma$ the subsets $\delta\cdot Z_{\gamma}$ and $Z_{\gamma}$ are disjoint: otherwise there are $x,y\in Z_{\gamma}$ with $y=\delta\cdot x$, and we have
\begin{equation*}
(\gamma\delta)\cdot x\ =\ \gamma\cdot y\ =\ \varphi(y)\ =\ \varphi(x)\ =\ \gamma\cdot x,
\end{equation*}
which is absurd since the action is free but $\gamma\delta\neq\gamma$.
Now recall that the measure $\mu$ is $2$-valued and $\varGamma$-invariant: it follows that each piece $Z_{\gamma}$, which has the same measure as $\delta\cdot Z_{\gamma}$ but is disjoint from it, is necessarily a $\mu$-nullset. We have thus that $(X,E)$ is covered by countably many nullsets, a contradiction.

(2.)
Let $f$ be an endomorphism of $(X,E_{\varDelta}^{(X,E)})$, whose graph is contained in the quotient relation $F/E_{\varDelta}^{(X,E)}\subseteq E_{\varGamma}^{(X,E)}/E_{\varDelta}^{(X,E)}$.
By the weak uniformization lemma \ref{lemma:WeakUniformization}, $f$ has a Borel lifting to an endomorphism $\varphi\in\End(X,E)$; let $\tilde{\varphi}$ and $Z_{\gamma}$, for $\gamma\in\varGamma$, be defined as in the proof of (1.): we will show that the unique $\gamma$ such that $\mu(Z_{\gamma})=1$ must belong to $N_{\varGamma}(\varDelta)$.
For every element $\delta$ of $\varDelta$, the intersection of the two full-measure subsets $Z_{\gamma}$ and $\delta\cdot Z_{\gamma}$ is nonempty, so there are $x,y\in Z_{\gamma}$ with $y=\delta\cdot x$. Since $x$ and $y$ are $E_{\varDelta}^{(X,E)}$-equivalent, their images by $\varphi$ have to be $E_{\varDelta}^{(X,E)}$-equivalent too, so, for some $\delta'\in\varDelta$,
\begin{equation*}
\gamma\delta\cdot x\;=\;\gamma\cdot y\;=\;\varphi(y)\;=\;\delta'\cdot\varphi(x)\;=\;\delta'\gamma\cdot x;
\end{equation*}
since the action of $\varGamma$ is free, necessarily $\gamma\delta=\delta'\gamma$. It follows that $\gamma\varDelta\gamma^{-1}\subseteq\varDelta$, so $\gamma\in N_{\varGamma}(\varDelta)$, as desired.
\qed

We now illustrate some easy applications of the preceding proposition, which are relevant for understanding the relationship of the various classes of countable equivalence relations discussed in \parref{sec:CountableEqRel} (see figure \ref{fig:ClassesCountableEqRelQuotientSpaces}).

\begin{ex}[A nonsmooth equivalence relation of index $2$]\label{ex:NonsmoothIndex2EqRel}
Let $E_0$ be the equivalence relation \eqref{eq:RelationE0} of eventual equality on the Cantor space \expset{\omega}{2};
it is well-known (see \cite{HKLGlimmEffrosDichotomy}, \cite{DJKStructureHyperfiniteEqRelations}) that $E_0$ is nonsmooth and that Haar measure $\mu_{(\expset{\omega}{2})}$, which is just the $\omega$-th power of the $\big(\frac{1}{2},\frac{1}{2}\big)$-measure on $2$, is $E_0$-nonatomic and ergodic ($0$-$1$ law): let $\mu=\mu_{(\expset{\omega}{2})}/E_0$ be the quotient measure.
Consider the componentwise action of the symmetric group on $2$ elements, $\Simm_2\cong\Z/2\Z$, on \expset{\omega}{2}:
\begin{equation}\label{eq:ActionSimm2Cantor}
\Simm_2\times\expset{\omega}{2}\to\expset{\omega}{2} \quad:\quad \sigma\cdot(x_n)\ =\ (\sigma x_n);
\end{equation}
this action preserves $E_0$-equivalence, so it induces an action of $\Simm_2$ on the quotient Borel space $(\expset{\omega}{2},E_0)$, which is easily seen to be Borel, free and $\mu$-preserving; its orbit equivalence is $E_0(\Simm_2)/E_0$, where $E_0(\Simm_2)$ is the equivalence relation generated by $E_0$ and $E^{\expset{\omega}{2}}_{\Simm_2}$.
By proposition \ref{prop:CriterionFreeMeasurePreservingActions}, $E_0(\Simm_2)/E_0$ is nonsmooth.
\end{ex}

\begin{ex}[The tail equivalence relation $E_t$ over $E_0$]\label{ex:EtOverE0}
Recall that the tail equivalence relation on the Cantor space is defined by
\begin{equation}\label{eq:RelationEt}
\seqn{x_n}\ E_t\ \seqn{y_n}\qquad\text{iff}\qquad\exists l,m\ \forall n\ (x_{l+n}=y_{m+n});
\end{equation}
in other words, $E_t$ is the equivalence generated by the \emph{shift endomorphism} $\sigma\in\End(\expset{\omega}{2})$, $\sigma(x_n)=(x_{n+1})$: $E_t=E(\sigma)$ (compare with \eqref{eq:TailEquivalenceEndomorphism}).
Note that the shift $\sigma$ preserves both $E_0$-equivalence and Haar measure $\mu_{(\expset{\omega}{2})}$, so it induces a Borel $\mu$-preserving endomorphism $\sigma'\in\End(\expset{\omega}{2},E_0)$.

We claim that $E_t/E_0$ is nonsmooth. Let $\expset{\omega}{2}_{*}$ be the $E_t$-invariant set of the aperiodic binary sequences, $\expset{\omega}{2}_{*}=\set{x\in\expset{\omega}{2}\ :\ \forall m,n\in\N\ (\sigma^m(x)=\sigma^n(x)\rightarrow m=n)}$: $\expset{\omega}{2}_{*}$ is a cocountable \BPi{0}{2} subset of \expset{\omega}{2}, so it is sufficient to study the restrictions of $E_0$ and $E_t$ to $\expset{\omega}{2}_{*}$. Observe now that the restriction of $\sigma'$ to $(\expset{\omega}{2}_{*},E_0)$ is an \emph{aperiodic automorphism}, so $E_t/E_0$ is the orbit equivalence of the Borel, free, $\mu$-preserving $\Z$-action induced by $\sigma'$: by proposition \ref{prop:CriterionFreeMeasurePreservingActions}, $E_t/E_0$ is nonsmooth.
\end{ex}

\begin{ex}[A smooth equivalence relation of index $3$ which is not enumerable]\label{ex:SmoothNonEnumerableIndex3EqRel}
Same notations as in example \ref{ex:NonsmoothIndex2EqRel}.
Consider the standard Borel space $X=\expset{\omega}{2}\times 2$ and the following two countable Borel equivalences on it:
\begin{align*}
&E\;=\;E_0(\Simm_2)\oplus E_0 &:&& &(x,i)E(y,j)\ \leftrightarrow\ \left\{
\begin{array}{l}
i=j=0\\
x E_0(\Simm_2) y
\end{array}
\right.\ \text{or}\ \left\{
\begin{array}{l}
i=j=1\\
x E_0 y
\end{array}
\right.\\
&F\;=\;E_0(\Simm_2)\times I(2) &:&& &(x,i)F(y,j)\ \leftrightarrow\ x E_0(\Simm_2) y.
\end{align*}
By construction, $F$ has index $3$ over $E$; moreover $F/E$ is smooth, since the set $(\expset{\omega}{2}\times\set{0})/E$ is a Borel transversal or, equivalently, the function
\begin{equation}
\tilde{\varphi}\;:\;X\to X\;:\;(x,i)\mapsto (x,0)
\end{equation}
induces a Borel selector $\varphi\in\End(X,E)$ for $F/E$ (see figure \ref{fig:SmoothNonEnumerableIndex3EqRel}).

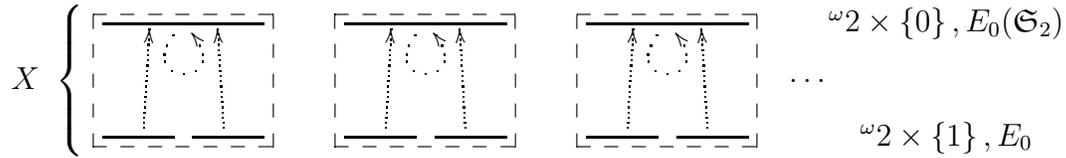
\begin{figure}[ht]
\begin{equation*}
X\ \left\{\ \vcenter{
\xymatrix@C=2em{
*+{}\ar@{-}[rr]_{}="ar1" & &*+{}&*+{}\ar@{-}[rr]_{}="ar2" & &*+{}&*+{}\ar@{-}[rr]_{}="ar3" & &*+{}&{\expset{\omega}{2}\times\set{0},E_0(\Simm_2)}\\
*+{}\ar@{-}[r]^{}="ar10" &\ar@{-}[r]^{}="ar11" &*+{}&*+{}\ar@{-}[r]^{}="ar20" &\ar@{-}[r]^{}="ar21" &*+{}&*+{}\ar@{-}[r]^{}="ar30" &\ar@{-}[r]^{}="ar31" &*+{}&{\expset{\omega}{2}\times\set{1},E_0}
\save "1,1"."2,3"*[F--]\frm{} \restore \save "1,4"."2,6"*[F--]\frm{} \restore
\save "1,7"."2,9"*[F--]\frm{}="last" \restore \save "last"!C+<4.5em,0pt>*{\ldots} \restore
\save \ar@{.>}"ar10";"ar1"+<-1.1em,-2pt> \ar@{.>}"ar11";"ar1"+<+1.1em,-2pt> \ar@{.>}@(dl,dr)"ar1";"ar1" \restore
\save \ar@{.>}"ar20";"ar2"+<-1.1em,-2pt> \ar@{.>}"ar21";"ar2"+<+1.1em,-2pt> \ar@{.>}@(dl,dr)"ar2";"ar2" \restore
\save \ar@{.>}"ar30";"ar3"+<-1.1em,-2pt> \ar@{.>}"ar31";"ar3"+<+1.1em,-2pt> \ar@{.>}@(dl,dr)"ar3";"ar3" \restore
}}\right.
\end{equation*}
\caption{Equivalence relations in example \ref{ex:SmoothNonEnumerableIndex3EqRel}: every $F$-class (dashed rectangle) splits into three $E$-classes (lines). The dotted arrows represent the selector $\varphi$.}
\label{fig:SmoothNonEnumerableIndex3EqRel}
\vspace{0.5cm}
\end{figure}

Suppose $F/E$ is enumerable. Consider the canonical projection $\pi:\frac{\expset{\omega}{2}}{E_0}\to\frac{\expset{\omega}{2}}{E_0(\Simm_2)}$ and the two Borel injections
\begin{align*}
&j_0\;:\;(\expset{\omega}{2},E_0(\Simm_2))\to(X,E) &\text{induced by} && &\tilde{\jmath}_0\;:\;x\mapsto (x,0),\\
&j_1\;:\;(\expset{\omega}{2},E_0)\to(X,E) &\text{induced by} && &\tilde{\jmath}_1\;:\;x\mapsto (x,1):
\end{align*}
their product $j_0\times j_1$ is a Borel embedding of $(\expset{\omega}{2},E_0(\Simm_2))\times(\expset{\omega}{2},E_0)$ into $(X,E)^2$, which maps the inverse relation $\pi^{-1}$ onto a subset $R$ of $F/E$. By the weak uniformization lemma \ref{lemma:WeakUniformization}, the enumerability assumption on $F/E$ insures the existence of a Borel uniformization $\psi_R$ of $R$: however, this means that $\psi=(j_0\times j_1)^{-1}[\psi_R]$ is a Borel right inverse of $\pi$, hence $\psi\circ\pi$ is a Borel selector for $E_0(\Simm_2)$ over $E_0$, contradicting the fact that $E_0(\Simm_2)/E_0$ is nonsmooth (example \ref{ex:NonsmoothIndex2EqRel}).
\end{ex}

\begin{remark}
By the representation theorem \ref{theo:FeldmanMooreQuotientSpaces}, a Borel equivalence relation is enumerable if and only if it is generated by countably many Borel \emph{automorphisms}. Generation by $\aleph_0$ countable-to-$1$ endomorphisms, on the other hand, is a much weaker condition for equivalences on quotient spaces, since Lusin\textendash Novikov uniformization no longer holds: the equivalence relation $F$ in example \ref{ex:SmoothNonEnumerableIndex3EqRel} is generated by a single Borel $3$-to-$1$ function $\varphi$ (a selector for $F$), and there is simply no countable family of Borel maps which can bring each point back in a uniform way to all its possible preimages by $\varphi$.
\end{remark}

\begin{ex}[An index-$3$ equivalence relation which is not countably generated]\label{ex:NotCtblyGeneratedIndex3EqRel}
(Louveau)
Consider the space \expset{\omega}{3} of infinite ternary sequences; following the discussion in example \ref{ex:NonsmoothIndex2EqRel} we define the equivalence relation $E_0$ of eventual equality and, for each subgroup $\Delta\leq\Simm_3$ of the symmetric group on $3$ elements, we consider the componentwise action of $\Delta$ on \expset{\omega}{3}, giving rise to the equivalence $E_0(\Delta)$ over $E_0$. The quotient equivalence relation $E_0(\Delta)/E_0$ is simply the orbit equivalence $E_{\Delta}^{(\expset{\omega}{3},E_0)}$ of the induced Borel action of $\Delta$ on $(\expset{\omega}{3},E_0)$.
Let $\mu_{(\expset{\omega}{3})}$ be Haar measure, i.e. the $\omega$-th power of the $\big(\frac{1}{3},\frac{1}{3},\frac{1}{3}\big)$-measure on $3$, and let $\mu=\mu_{(\expset{\omega}{3})}/E_0$ be the quotient measure; $\mu$ is $2$-valued and $\Simm_3$-invariant.
In order to apply proposition \ref{prop:CriterionFreeMeasurePreservingActions} we consider the restrictions of the previous equivalence relations to the cocountable $E_0(\Simm_3)$-invariant \BPi{0}{2} subset $\expset{\omega}{3}_{*}\subseteq\expset{\omega}{3}$ of aperiodic sequences: in fact, the action of $\Simm_3$ on $(\expset{\omega}{3}_{*},E_0)$ is now \emph{free}. Let $\Simm_2$ be identified by the subgroup of $\Simm_3$ generated by a transposition: $[\Simm_3:\Simm_2]=3$ and the normalizer $N_{\Simm_3}(\Simm_2)$ of $\Simm_2$ is $\Simm_2$ itself. It follows that $E_0(\Simm_3)$ has index $3$ over $E_0(\Simm_2)$, but any Borel countably generated subequivalence relation of $E_0(\Simm_3)/E_0(\Simm_2)$ has to be almost equal to the identity (proposition \ref{prop:CriterionFreeMeasurePreservingActions}), so $E_0(\Simm_3)$ is not countably generated over $E_0(\Simm_2)$.
\end{ex}

\begin{ex}[An index-$3$ equivalence relation which cannot be generated by countably many involutions]\label{ex:NotInvolutionGeneratedIndex3EqRel}
Same notations as in example \ref{ex:NotCtblyGeneratedIndex3EqRel}.
Consider the quotient relation $E_0(\Z/3\Z)/E_0=E_{\Z/3\Z}^{(\expset{\omega}{3},E_0)}$, where $\Z/3\Z$ is identified with the subgroup of $\Simm_3$ generated by an element of order $3$ (a rotation); observe that the action of $\Z/3\Z$ on $(\expset{\omega}{3},E_0)$ is free, so pick the associated Borel cocycle $\theta$ (lemma \ref{lemma:ExistenceOfBorelCocyclesFreeActions}).
We will show that every involution $f\in\Aut(\expset{\omega}{3},E_0)$ with $f\subseteq E_0(\Z/3\Z)/E_0$ has $\mu$-null support, i.e. $f(x)=x$ for $\mu$-almost every $x$: this clearly implies that $E_0(\Z/3\Z)/E_0$ cannot be generated by countably many involutions.
We follow the same reasoning used in the proof of proposition \ref{prop:CriterionFreeMeasurePreservingActions}: given an endomorphism $f\in\End(X,E_0)$ with $f\subseteq E_0(\Z/3\Z)/E_0$, let $\tilde{f}(x)=\theta(x,f(x))$ and, for $i\in\Z/3\Z$, let $Z_i=\set{x\in (\expset{\omega}{3},E_0):\tilde{f}(x)=i}$; then $\mu(Z_i)=1$ for exactly one $i$ and $\mu(Z_j)=0$ for the other $j\neq i$.
Suppose $f^2=\id$: then $f[Z_0]=Z_0$, $f[Z_1]=Z_2$ and $f[Z_2]=Z_1$; since, by construction, $f[Z_i]=i\cdot Z_i$ and $\mu$ is $\Z/3\Z$-invariant, the only possibility is that $\mu(Z_0)=1$: this leads to the desired conclusion, since $f$ is the identity on $Z_0$.
\end{ex}

\commentout{
\begin{ex}[A finite orbit equivalence relation which cannot be generated by countably many automorphisms of finite order]\label{ex:NoFiniteOrderGeneratorsFiniteEqRel}

\end{ex}
}


\bibliographystyle{amsalpha}
\bibliography{../Main}

\providecommand{\bysame}{\leavevmode\hbox to3em{\hrulefill}\thinspace}
\providecommand{\MR}{\relax\ifhmode\unskip\space\fi MR }
\providecommand{\MRhref}[2]{%
  \href{http://www.ams.org/mathscinet-getitem?mr=#1}{#2}
}
\providecommand{\href}[2]{#2}
\begin{thebibliography}{HKL90}

\bibitem[DJK94]{DJKStructureHyperfiniteEqRelations}
R.~Dougherty, S.~Jackson, and A.~S. Kechris, \emph{The structure of hyperfinite
  {B}orel equivalence relations}, Trans. Amer. Math. Soc. \textbf{341} (1994),
  no.~1, 193--225.

\bibitem[Eff65]{EffTransformationGroupsCAlgebras}
Edward~G. Effros, \emph{Transformation groups and {$C\sp{\ast} $}-algebras},
  Ann. of Math. (2) \textbf{81} (1965), 38--55.

\bibitem[Eff81]{EffPolishTransfGroupsClassification}
\bysame, \emph{Polish transformation groups and classification problems},
  General topology and modern analysis (Proc. Conf., Univ. California,
  Riverside, Calif., 1980), Academic Press, New York, 1981, pp.~217--227.

\bibitem[HKL90]{HKLGlimmEffrosDichotomy}
L.~A. Harrington, A.~S. Kechris, and A.~Louveau, \emph{A {G}limm-{E}ffros
  dichotomy for {B}orel equivalence relations}, J. Amer. Math. Soc. \textbf{3}
  (1990), no.~4, 903--928.

\bibitem[JKL02]{JKLCountableBorelEqRelations}
S.~Jackson, A.~S. Kechris, and A.~Louveau, \emph{Countable {B}orel equivalence
  relations}, J. Math. Log. \textbf{2} (2002), no.~1, 1--80.

\bibitem[Kec95]{KecClassicalDST}
Alexander~S. Kechris, \emph{Classical descriptive set theory}, Graduate Texts
  in Mathematics, vol. 156, Springer-Verlag, New York, 1995.

\bibitem[KM04]{KMTopicsOrbitEquivalence}
Alexander~S. Kechris and Benjamin~D. Miller, \emph{Topics in orbit
  equivalence}, Lecture Notes in Mathematics, vol. 1852, Springer-Verlag,
  Berlin, 2004.

\bibitem[MK80]{MKInfiniteGamesEffectiveDST}
D.~A. Martin and A.~S. Kechris, \emph{Analytic sets}, ch.~Infinite games and
  effective descriptive set theory, pp.~403--470, Academic Press Inc. [Harcourt
  Brace Jovanovich Publishers], 1980, Lectures delivered at a Conference held
  at University College, University of London, London, July 16--29, 1978.

\bibitem[Sil80]{SilCountingNumberEquivalenceClasses}
Jack~H. Silver, \emph{Counting the number of equivalence classes of {B}orel and
  coanalytic equivalence relations}, Ann. Math. Logic \textbf{18} (1980),
  no.~1, 1--28.

\bibitem[Wei84]{WeiMeasurableDynamics}
Benjamin Weiss, \emph{Measurable dynamics}, Conference in modern analysis and
  probability (New Haven, Conn., 1982), Contemp. Math., vol.~26, Amer. Math.
  Soc., Providence, RI, 1984, pp.~395--421.

\end{thebibliography}

\end{document}